\documentclass[12pt]{article}
\usepackage[numbers]{natbib}
\usepackage[margin=1in]{geometry}
\usepackage{amssymb}
\usepackage{amsmath}
\usepackage{graphicx,tikz,bookmark,epsfig,rotating,colortbl}%
\usepackage{multirow, anyfontsize}%

\usepackage{amsthm}%
\usepackage{mathrsfs}%
\usepackage[title]{appendix}%
\usepackage{xcolor}%
\usepackage{textcomp}%
\usepackage{manyfoot}%
\usepackage{booktabs}%
\usepackage{authblk}
\usepackage{listings,bm}
\usepackage{comment}
\usepackage{float}
\usepackage{hhline}
\def\tsc#1{\csdef{#1}{\textsc{\lowercase{#1}}\xspace}}
\tsc{WGM}
\tsc{QE}

\begin{document}

\let\WriteBookmarks\relax
\def\floatpagepagefraction{1}
\def\textpagefraction{.001}
\definecolor{awesome}{rgb}{1.0, 0.13, 0.32}
\definecolor{electricindigo}{rgb}{0.44, 0.0, 1.0}
\definecolor{burgundy}{rgb}{0.5, 0.0, 0.13}
\definecolor{amethyst}{rgb}{0.6, 0.4, 0.8}
\newtheorem{proposition}{Proposition}
\newtheorem{definition}{Definition}%
\newtheorem{corollary}{Corollary}[section]
\newtheorem{lemma}{Lemma}[section]
\newtheorem{theorem}{Theorem}
\newtheorem{condition}{Condition}



\title{\textbf {Infinitely Many Almost-Uniform $3$-Qualitative Independence Hypergraphs are Cores}  }



%



\author[1]{Raina Mary Thomas}
\author[2]{Yasmeen Akhtar}

\affil[1,2]{ Department of Mathematics, BITS Pilani K K Birla Goa Campus, Goa 403726, India}

\affil[1]{\texttt{p20220019@goa.bits-pilani.ac.in}}
\affil[2]{\texttt{yasmeena@goa.bits-pilani.ac.in}}

\date{}

\maketitle


\begin{abstract}
  Qualitative independence is a property associated with the column vectors of well-known combinatorial designs known as covering arrays and qualitative independence hypergraphs characterize the existence of covering arrays on hypergraphs. The core of a hypergraph $H$ is an interesting structural characteristic to study, as it is the unique (up to isomorphism), minimal subhypergraph $H^\bullet \subseteq H$ which is homomorphically equivalent to $H$. We use the vertex transitivity and strong independence number of almost-uniform $3$-qualitative independence hypergraphs, $3\text{-}AUQI(2n+1,2)$, to show that $3\text{-}AUQI(2n+1,2)$ are cores for infinitely many $n$. Consequently, we establish that $3\text{-}QI(9,2)$ is a core and also obtain a lower bound on the strong chromatic number of $3$-qualitative independence hypergraphs $3\text{-}QI(2n+1,2)$. 
\end{abstract}







\section{Background}\label{sec:1}
For positive integers $n$, $t$ such that $t \geq 2$, $n \geq g^t$, and a set $S$ of $g$ symbols, a set of strings of length $n$ and entries from $S$ is said to be \emph{$t$-qualitatively independent} ($t$-QI) if for any $t$-subset, $\{v_{1}, v_{2}, \ldots , v_{t}\}$, of strings and any ordered $t$-tuple say $(g_1, g_2, \ldots, g_t)\in S^t$, there exists an index $j$ such that $v_{i}[j] = g_i$ for $i=1, 2, \ldots, t$. Each of these strings $v$ can be represented in terms of partitions, by defining a corresponding partition of $[n]=\{1,2,\ldots, n\}$ to be $R = (M_{1}, M_{2}, \ldots, M_{g})$, where $i \in M_{k}$ whenever $v[i]=k-1$, ($1 \leq i \leq n, 1 \leq k \leq g$). For an $n$-element set, $N$, a family of partitions of $N$ into $g$ classes, $\mathcal{P} = \{R_1, R_2, \ldots, R_k\}$, where $R_j = (M^j_1, M^j_2, \ldots, M^j_g)$ for $j = 1, 2, \ldots, k$, is said to be \emph{$t$-qualitatively independent} if, for every $t$ of the partitions $R_{i_1}, R_{i_2}, \ldots, R_{i_t}$, and for every choice of classes $M^{i_j}_s$ (with $1 \leq j \leq t, 1 \leq s \leq g$), $\displaystyle \bigcap_{j=1}^{t}M^{i_j}_s \neq \emptyset$, for all choices of $s$ and distinct choices of $i_j$. Thus, a set of partitions is $t$-qualitatively independent if and only if the corresponding strings are $t$-qualitatively independent. Figure \ref{QI} shows a set of $3$-QI binary strings of length $8$ along with their corresponding $3$-QI $2$-partitions of $[8]$. The significance of this notion in search theory lies in the consequence that, after knowing if an unknown $x$ is in a set $A$ or not, we cannot decide the same question for another set $B$, \emph{independently}, based on the membership of $x$ in $A$; hence justifying the nomenclature \cite{katona2004strongqi}. \\

    Motivated by a problem in \emph{probability theory} to design a series of experiments which are more or less independent, R\'{e}nyi \cite{renyi} posed the problem to determine the maximum size of a family of $t$-QI partitions for given $g$ and $n$, generalizing an early investigation done by Marczewski \cite{marczewski1948independance} for the case $g=2$. R\'{e}nyi proved that it is $\dbinom{n-1}{ \frac{n}{2} -1}$, for even $n$, using the well-known Sperner's theorem \cite{sperner1928satz}. 
    Later, Kleitman and Spencer \cite{kleitman1973families}, Katona \cite{katona1973two}, Brace and Daykin \cite{brace1972sperner} and Bollob\'as \cite{bollobas1973sperner} independently determined that it is $\dbinom{n-1}{\lfloor \frac{n}{2}\rfloor -1}$ for all $n$, when $t = g = 2$. They proved the case for odd $n$, by additionally using the Erd\H{o}s-Ko-Rado Theorem \cite{erdos1961intersection}, an important result in \emph{extremal set theory}. 
    Poljak et al. \cite{poljak1983qualitatively} gave a constructive procedure to find the maximum number of $t$-QI subsets of an $n$-set for $t>2$, in contrast to the existential proofs given by Kleitman, Spencer and Katona. Shortly after this, Alon \cite{alon1986explicit} used the theory of error correcting codes to give an explicit construction of an exponential-sized family of subsets of an $n$-set that is $t$-QI. The size of the constructed family was at least $2^{cn}$, where $c$ is a positive constant independent of $n$. Poljak and Tuza \cite{poljak1989maximum} constructed larger $2$-QI collections, improved the order of magnitude of the previously known upper bounds on the maximum number of $t$-QI partitions, and gave a general lower bound using probabilistic arguments. When $g=2$, $2$-QI partitions are also studied under the name \emph{crossing} sets. Two subsets $A$ and $B$ of a ground set $X$ are \emph{crossing} if none of the four sets $A \backslash B$, $B \backslash A$, $A \cap B$, $X \backslash (A \cup B)$ are empty, equivalently, $A$ and $B$ are crossing if their Venn diagram has no empty regions \cite{karzanov1978flow, lomonosov1985combinatorial, brass2005research, dress2012basic, tomon2026cross}. Some other concepts related to qualitative independence found in literature are separating partitions \cite{poljak1983qualitatively}, weak qualitative dependence \cite{korner2002delta} and strong qualitative independence \cite{katona2004strongqi}.\\

    \begin{figure}\begin{center}
    \begin{tabular}{|c|c|}
 \hline
         \textcolor{blue}{0000}\textcolor{awesome}{1111}&\{\textcolor{blue}{\{1, 2, 3, 4\}}, \textcolor{awesome}{\{5, 6, 7, 8\}}\}\\

          \hline
         00110011~&\{\{1, 2, 5, 6\}, \{3, 4, 7, 8\}\} \\
         \hline
          01010101~&\{\{1, 3, 5, 7\}, \{2, 4, 6, 8\}\} \\
         \hline
          01101001~&\{\{1, 4, 6, 7\}, \{2, 3, 5, 8\}\} \\
         \hline
      
    \end{tabular}\end{center}\caption{A set of $3$-QI strings of length $8$ and corresponding $3$-QI partitions of $[8]$}
    \label{QI}\end{figure}

\begin{figure}
    \centering
    \begin{tikzpicture}[scale=1]
        \node[draw, black, circle,thick,minimum size=0.2cm,inner sep=0pt, fill=black] at (0,0) {};
        \node[draw, black, circle,thick,minimum size=0.2cm,inner sep=0pt, fill=black] at (0,2) {};
        \node[draw, black, circle,thick,minimum size=0.2cm,inner sep=0pt, fill=black] at (2,0) {};
        \node[draw, black, circle,thick,minimum size=0.2cm,inner sep=0pt, fill=black] at (2,2) {};
        \node[draw, black, circle,thick,minimum size=0.2cm,inner sep=0pt, fill=black] at (2.5,1) {};
        \node[draw, black, circle,thick,minimum size=0.2cm,inner sep=0pt, fill=black] at (-0.5,1) {};
        \node[draw, black, circle,thick,minimum size=0.2cm,inner sep=0pt, fill=black] at (1,2.5) {};
        \node[draw, black, circle,thick,minimum size=0.2cm,inner sep=0pt, fill=black] at (1,-0.5) {};
        \node[draw, black, circle,thick,minimum size=0.2cm,inner sep=0pt, fill=black] at (4,1.5) {};
        \node[draw, black, circle,thick,minimum size=0.2cm,inner sep=0pt, fill=black] at (0.5,-2) {};
        \draw (0, 0)--(0, 2) -- (2, 2) --(2, 0)--(0, 0);
        \draw (2.5,1)--(1,2.5)--(-0.5,1)--(1, -0.5)--(2.5,1);
        \draw (-0.5,1)--(2.5,1)--(4,1.5)--(2,0)--(0.5,-2)--(1,-0.5)--(1,2.5);
        \draw (-0.5,1)--(0,0)--(1,-0.5)--(2,0)--(2.5,1)--(2,2)--(1,2.5)--(0,2)--(-0.5,1);
        \draw (2,2)--(4,1.5);
        \draw (0,0)--(0.5,-2);
        \draw (0, 2) --(4,1.5)--(-0.5, 1);
        \draw (0, 0)--(4, 1.5);
        \draw (1, 2.5)--(0.5,-2)--(2,2);
        \draw (0.5,-2)--(0,2);
        \node at (1, 3) {00111};
        \node at (2.75, 2) {01001};
        \node at (2.75, 1.25) {01110};
        \node at (2.75, 0){00011};
        \node at (4.5, 1) {00101};
        \node at (1.75, -0.5) {01101};
        \node at (1.25, -2) {01010};
        \node at (-0.75, 0) {00110};
        \node at (-1.25, 1) {01011};
        \node at (-0.75, 2){01100};
        
    \end{tikzpicture}
   
    \caption{$2\text{-}QI(5,2)\cong 2\text{-}AUQI(5,2)$}
    \label{fig:qi}
\end{figure}
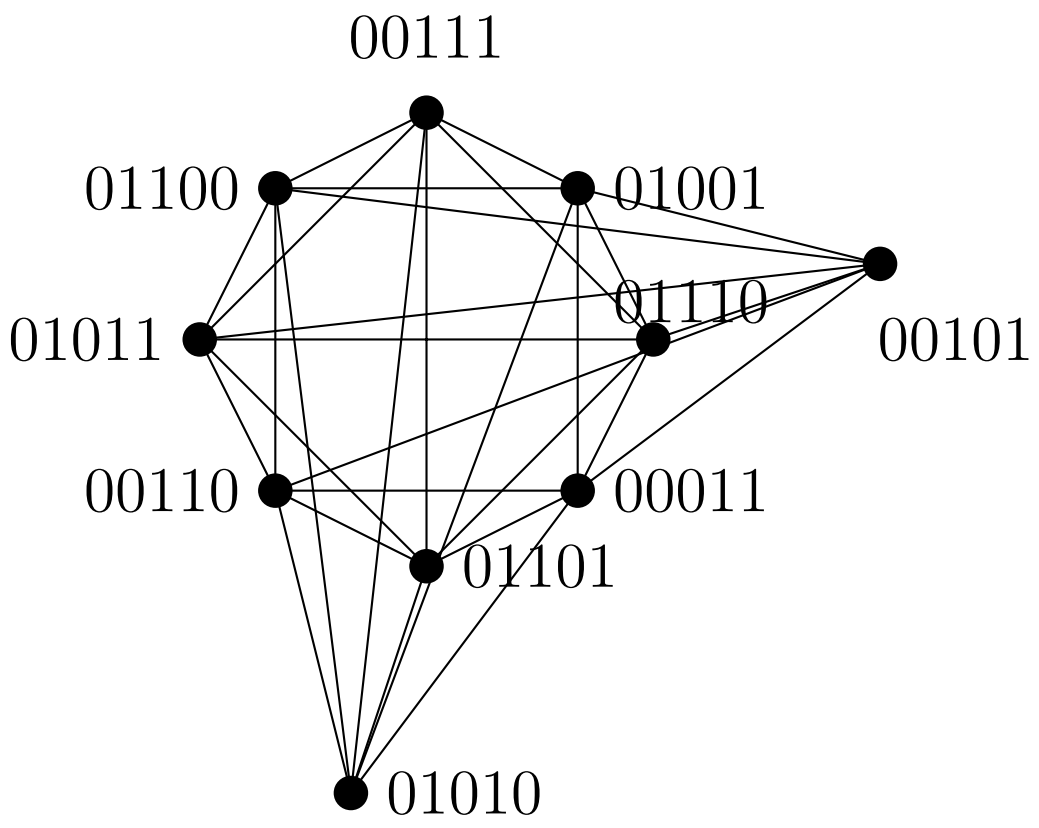

The qualitative independence property forms the basis of \emph{qualitative independence (QI) hypergraphs}, which are studied in the context of \emph{design theory} to characterize the existence of covering arrays on hypergraphs \cite{meagher2005covering,Raaphorst2013, raaphorst2018variable}. A \emph{hypergraph $H$} is a pair $H = (V, E)$, where $V = \{v_1, v_2, \ldots , v_k\}$ is a set of vertices, and $E = \{e_1, e_2, \ldots , e_m\}$ is a set of non-empty subsets of $V$, called \emph{hyperedges}. A \emph{$t$-uniform} hypergraph is one in which all the hyperedges have the same cardinality $t$. Figure \ref{example} shows a $3$-uniform hypergraph $H$ with vertex set $\{a,b,c,d,e,f,g\}$ and hyperedge set $\{\{a,b,c\}, \{b,c,e\}, \{b,c,d\}, \{b,c,f\}, \{d,e,f\} \{e,f,g\}\}$. A \emph{$t$-qualitative independence hypergraph}, $t\text{-}QI(n, g) = (V, E)$, is a $t$-uniform hypergraph, where $V$ is the set of all $g$-ary strings of length $n$ with entries from $S$ and with the property that every symbol appears at least $g^{t-1}$ times and the first appearance of every symbol is in the lexicographic order; a set of $t$ vertices forms a hyperedge if they are $t$-qualitatively independent. Equivalently, the vertices may be represented as partitions of an $n$-set into $g$ classes such that the indices where a particular symbol occurs in the string belong to the same class in the partition. Figure \ref{fig:qi} shows $2\text{-}QI(5,2)$.\\
\begin{figure}
    \centering
    \begin{tikzpicture}[scale=4]

        \draw[
            rotate around={0:(0,0)},
            black,
            fill=cyan!40,
            fill opacity=0.4
        ]
        (0.05,0.8) ellipse (0.12cm and 0.50cm);

        \draw[
            rotate around={0:(0,0)},
            black,
            fill=red,
            fill opacity=0.55
        ]
        (0.05,0.4) ellipse (0.15cm and 0.48cm);

        \draw[
            black,
            rounded corners=13mm,
            fill=green!60,
            fill opacity=0.45
        ]
        (0.1,1.05)--(0.1,0.05)--(-0.6,-0.2)--cycle;
                \draw[
            black,
            rounded corners=13mm,
            fill=orange!60,
            fill opacity=0.35
        ]
        (0.45,0) ellipse (0.5cm and 0.15cm);

        \draw[
            black,
            rounded corners=13mm,
            fill=electricindigo!60,
            fill opacity=0.6
        ]
        (0,1.05)--(0,0.05)--(0.6,-0.2)--cycle;

        \draw[
            black,
            rounded corners=13mm,
            fill=gray!60,
            fill opacity=0.6
        ]
       (0,0) ellipse (0.5cm and 0.15cm);

        \node[
            draw, black, circle, thick,
            minimum size=0.1cm,
            inner sep=0pt,
            fill=black
        ] at (0.05,0.05) {};

        \node[
            draw, black, circle, thick,
            minimum size=0.1cm,
            inner sep=0pt,
            fill=black
        ] at (0.05,0.42) {};

        \node[
            draw, black, circle, thick,
            minimum size=0.1cm,
            inner sep=0pt,
            fill=black
        ] at (0.05,0.80) {};

        \node[
            draw, black, circle, thick,
            minimum size=0.1cm,
            inner sep=0pt,
            fill=black
        ] at (0.05,1.18) {};

        \node[
            draw, black, circle, thick,
            minimum size=0.1cm,
            inner sep=0pt,
            fill=black
        ] at (-0.28,0.05) {};

        \node[
            draw, black, circle, thick,
            minimum size=0.1cm,
            inner sep=0pt,
            fill=black
        ] at (0.4,0.05) {};
        \node[
            draw, black, circle, thick,
            minimum size=0.1cm,
            inner sep=0pt,
            fill=black
        ] at (0.69,0.05) {};

        \node at (0.2,1.18) {\large $a$};
        \node at (0.2,0.80) {\large $b$};
        \node at (0.17,0.42) {\large $c$};

        \node at (0.05,0) {\large $e$};
        \node at (-0.32,0) {\large $d$};
        \node at (0.37,-0.02) {\large $f$};
        \node at (0.69,-0.02) {\large $g$};

        \node at (0.05,-0.42) {$\boldsymbol{H}$};

    \end{tikzpicture}
    \begin{tikzpicture}[scale=4]

        \draw[
            rotate around={0:(0,0)},
            black,
            fill=red,
            fill opacity=0.55
        ]
        (0.05,0.4) ellipse (0.15cm and 0.48cm);

        \draw[
            black,
            rounded corners=13mm,
                        fill=green!60,
            fill opacity=0.45
        ]
        (0.1,1.05)--(0.1,0.05)--(-0.6,-0.2)--cycle;

        \draw[
            black,
            rounded corners=13mm,
            fill=electricindigo!60,
            fill opacity=0.6
        ]
        (0,1.05)--(0,0.05)--(0.6,-0.2)--cycle;

        \draw[
            black,
            rounded corners=13mm,
            fill=gray!60,
            fill opacity=0.6
        ]
       (0,0) ellipse (0.5cm and 0.15cm);

        \node[
            draw, black, circle, thick,
            minimum size=0.1cm,
            inner sep=0pt,
            fill=black
        ] at (0.05,0.05) {};

        \node[
            draw, black, circle, thick,
            minimum size=0.1cm,
            inner sep=0pt,
            fill=black
        ] at (0.05,0.42) {};

        \node[
            draw, black, circle, thick,
            minimum size=0.1cm,
            inner sep=0pt,
            fill=black
        ] at (0.05,0.80) {};

        \node[
            draw, black, circle, thick,
            minimum size=0.1cm,
            inner sep=0pt,
            fill=black
        ] at (-0.28,0.05) {};

        \node[
            draw, black, circle, thick,
            minimum size=0.1cm,
            inner sep=0pt,
            fill=black
        ] at (0.4,0.05) {};

        
        \node at (0.2,0.80) {\large $b$};
        \node at (0.16,0.42) {\large $c$};

        \node at (0.05,0) {\large $e$};
        \node at (-0.32,0) {\large $d$};
        \node at (0.37,-0.02) {\large $f$};

        \node at (0.05,-0.42) {$\boldsymbol{H^\bullet}$};
\node at (1.2, 0.40) {\begin{tabular}{>{\columncolor[gray]{0.90}}c>{\columncolor[gray]{0.90}}c>{\columncolor[gray]{0.90}}c>{\columncolor[gray]{0.90}}c>{\columncolor[gray]{0.90}}c}
					 \rowcolor{amethyst!30!}   \textbf{$b$}& \textbf{$c$} &  \textbf{$d$} & \textbf{$e$} & \textbf{$f$}\\		
							0 & 0 & 0 & 0 & 0\\
							0 & 1 & 0 & 0 & 1\\
							0 & 0 & 1 & 1 & 1\\
							0 & 1 & 1 & 1 & 0\\
							1 & 0 & 0 & 1 & 1\\
							1 & 1 & 0 & 1 & 0\\
							1 & 0 & 1 & 0 & 0\\
							1 & 1 & 1 & 0 & 1\\
							
					\end{tabular}};
    \end{tikzpicture}
    
    \caption{A hypergraph $H$, its core $H^\bullet$, and a $CA(8, H^\bullet, 2)$}
    \label{example}
\end{figure}

A \emph{covering array on a hypergraph} $H = (V, E)$, of size $n$ and $g$ symbols, denoted $CA(n, H, g)$, is an $n \times |V|$ array with entries from $S$ and columns indexed by the vertices in $V$ such that: for any hyperedge $e \in E$ with $|e| = m$, the columns in the sub-array corresponding to the vertices in $e$ are $m$-qualitatively independent. Therefore, covering arrays (CAs) are also known as \emph{qualitatively independent families}. They are used in combinatorial testing as a cost-effective and efficient tool to design test suites in various domains \cite{appl1994compressing,nie2011survey,applstevens1998efficient,drug,survey2019methods}. Figure \ref{example} shows a hypergraph $H^\bullet$ with vertex set $\{b,c,d,e,f\}$, hyperedge set $\{\{b,c,d\},\{b,c,e\},\{b,c,f\}, \{d,e,f\}\}$ and a covering array on it, $CA(8, H^\bullet, 2)$. Raaphorst et al. \cite{raaphorst2018variable} showed that for a $t$-uniform hypergraph $H$ and positive integers $g$ and $n$, a covering array $CA(n, H, g)$ exists if and only if there is a hypergraph homomorphism from $H$ to $t\text{-}QI(n, g)$, extending the work done by Meagher and Stevens \cite{meagher2005covering} who introduced and characterized the existence of covering arrays on graphs. The size of optimal CAs on QI graphs and some QI hypergraphs with binary symbols are given in \cite{meagher2005covering} and \cite{raina} respectively.\\

A \emph{uniform} $t$-qualitative independence hypergraph, $t\text{-}UQI(n, g)$, is an induced subhypergraph of $t\text{-}QI(n,g)$ whose vertices have the property that every symbol appears equally many times in the string. Cliques in $t\text{-}UQI(n, g)$ correspond to \emph{balanced} CAs, that is, a CA with the property that each symbol occurs the same number of times in each column. Raaphorst et al. \cite{raaphorst2018variable} proved that a balanced CA on a hypergraph $H$ with size $n$ and $g$ symbols exists if and only if there exists a hypergraph homomorphism from $H$ to $t\text{-}UQI(n, g)$. An \emph{almost-uniform} $t$-qualitative independence hypergraph, $t\text{-}AUQI(n, g)$, is an induced subhypergraph of $t\text{-}QI(n,g)$ whose vertices have the property that if $n = gs+r$ for some integers $s$ and $0 < r < g$, then $(g-r)$ of the symbols appear $s$ times and the remaining $r$ appears $(s+1)$ times in the string. Figure \ref{fig:qi} shows the smallest almost-uniform QI hypergraph, namely the graph $2\text{-}AUQI(5,2)\cong 2\text{-}QI(5,2)$. An \emph{almost-balanced} CA on a hypergraph $H$ with size $n$ and $g$ symbols exists if and only if there exists a hypergraph homomorphism from $H$ to $t\text{-}AUQI(n, g)$. These existence characterizations motivate the study of structural aspects of QI hypergraphs, in particular, their \emph{cores}. \\
  
 The core of a hypergraph $H$ is an interesting structural characteristic to study, as it is the minimal subhypergraph $H^\bullet \subseteq H$ which is homomorphically equivalent to $H$. It preserves the size of the maximum $t$-cliques of the hypergraph as well as its strong chromatic number. Every finite hypergraph has a unique core up to isomorphism \cite{micciancio2004using}. Cores of graphs have been extensively studied from the 1960s under various names such as minimal graphs, reduced graphs, automorphic graphs, retract-rigid graphs, etc., in diverse contexts \cite{bauslaugh1995core,fellner1982minimal, hell1992core,  nevsetvril1978classescore, royle2001algebraic}. Versions of cores with regard to relational structures have been found to be useful in database theory for conjunctive-query processing and in data exchange problems to find the best universal solution \cite{chandra1977optimal, fagin2003data}. The decision problem of CORE IDENTIFICATION, that is, ``given graphs $G_1$ and $G_2$, is $G_2$ the core of $G_1$?" is DP-complete \cite{fagin2003data} and the decision problem of CORE RECOGNITION, that is, ``given a graph $G$, is $G$ a core?" is co-NP-complete \cite{hell1992core}. Some well-known examples of cores include complete graphs $K_n$, odd cycles, odd wheels, Kneser graphs (including the Petersen graph), $k$-critical graphs, and the icosahedron, etc., whereas even cycles and complete graphs with one edge removed are not cores. Johnson graphs are either cores or have a complete core \cite{godsil2011cores}. For all positive integers $k$, $k$-uniform complete hypergraphs are cores. Cameron and Kazanidis \cite{cameron2008cores} discuss the cores of nonedge-transitive graphs.\\

From Raaphorst et al.'s characterization \cite{raaphorst2018variable}, if the core of $t\text{-}QI(n, g)$ is $t\text{-}UQI(n, g)$ for even $n$ and $t\text{-}AUQI(n, g)$ for odd $n$, then the corresponding existence of homomorphisms are sufficient to conclude that the existence of CAs with given parameters ensures the existence of, respectively, balanced and almost-balanced CAs with the same parameters. This motivates the question: \emph{Are $t\text{-}UQI(n, g)$ and $t\text{-}AUQI(n, g)$ cores?} Meagher and Stevens \cite{meagher2005covering} showed that the core of $2\text{-}QI(n, 2)$ is $2\text{-}UQI(n, 2)$, when $n$ is even, and $2\text{-}AUQI(n, 2)$, when $n$ is odd. In particular, $2\text{-}UQI(n, 2)$ and $2\text{-}AUQI(n, 2)$ are cores. Godsil and Newman \cite{godsilnewman33} showed that $2\text{-}QI(9,3)$ is a core. Thomas and Akhtar \cite{raina2026coresuniformalmostuniform3qualitative} showed that $3\text{-}QI(8,2)$ is a core and established a structural connection between QI hypergraphs and merged Johnson graphs. They also gave sufficient conditions under which $3\text{-}UQI(n,2)$, for even $n\geq 10$ and $3\text{-}AUQI(n, 2)$ for odd $n\geq 9$ are cores. To the best of our knowledge, the cores of $t\text{-}QI(n, 2)$ for $t>2$, except for $n=8$ and $2\text{-}QI(n, g)$ for $g\geq 3$ except when $n=9$, have not yet been identified. In this paper, we establish that $3\text{-}QI(9,2)$ is a core and prove that $3\text{-}AUQI(n, 2)$ is a core for infinitely many odd $n$, specifically whenever $n=p-2$, for prime $p\geq 11$. Equivalently, we conclude that the sufficient condition (based on merged Johnson graphs) given in \cite{raina2026coresuniformalmostuniform3qualitative} under which $3\text{-}AUQI(n, 2)$ is a core, holds true for such $n$.\\

The remainder of the article is organized as follows. Section \ref{sec:2} discusses some properties exhibited by vertex-transitive hypergraphs. In Section \ref{sec:3}, we determine the strong independence number of $3\text{-}AUQI(2n+1,2)$ for $n\geq 4$ and apply this to give a lower bound on the strong chromatic number of $3\text{-}QI(2n+1,2)$. We then use these results to establish that infinitely many $3$-AUQI hypergraphs and the associated merged Johnson graphs are cores in Section \ref{sec:4}; and conclude with some comments on open cases and future direction in Section \ref{sec:5}.

\section{Preliminaries}\label{sec:2}

In this section, we give some basic definitions and preliminary details on cores and vertex-transitive hypergraphs. A \emph{subhypergraph} $H'=(V', E')$ of a hypergraph $H=(V,E)$ is a hypergraph such that $V'\subseteq V$, $E' \subseteq E \cap \mathcal{P}(V')$, where $\mathcal{P}(V')$ is the power set of $V'$. If $E' = E \cap \mathcal{P}(V')$, then $H'$ is said to be an \emph{induced} subhypergraph. In Figure \ref{example}, $H^\bullet$ is an induced subhypergraph of $H$. A \emph{$t$-complete hypergraph of order $n$},  denoted $K_n^t$, is a $t$-uniform hypergraph $H=(V,E)$ such that $|V| = n$ and  $E$ consists of all the $t$-subsets of $V$. Given a hypergraph $H$, a \emph{$t$-clique} of size $m$ in $H$ is a subset $C \subseteq V$ of size $m$ such that the induced subhypergraph on $C$ is $K_m^t$. We write $\omega_t(H)$ to denote the size of the largest $t$-clique in $H$. For example, $\omega_3(H)=3$ for $H$ in Figure \ref{example}. The \emph{$k$-section} of a hypergraph $H=(V, E)$ is defined to be a hypergraph $[H]_k$ with the same vertex set and whose hyperedges are the sets $F \subseteq V$ satisfying either $|F|=k$ and $F \subset e$ for some $e \in E$; or $|F|<k$ and $F=e$ for some $e \in E$. The \emph{strong chromatic number} of a hypergraph $H$, denoted $\chi_S(H)$, is the minimum number of colours required to colour the vertices of $H$ such that no two vertex in a hyperedge are assigned the same colour. When $H$ is a graph, this corresponds to the usual \emph{chromatic number}. The vertices in the hypergraph $H$ in Figure \ref{example} can be coloured by assigning the same colour to $a$, $d$ and $g$, and distinct colours to the rest of the vertices, leading to $\chi_S(H)=5$.\\

For two hypergraphs $H_1=(V_1,E_1)$ and $H_2=(V_2, E_2)$, a \emph{hypergraph homomorphism} from $H_1$ to $H_2$, denoted $f : H_1 \rightarrow H_2$, is a map $f : V_1 \rightarrow V_2$ such that for every hyperedge $e = \{v_1,\ldots, v_m\} \in E_1$, there exists a hyperedge $e'\in E_2$ such that  the image $f(e) = \{f(v_1), \ldots,f(v_m)\} \subseteq e'$ and $|e|$ = $|f(e)|$. Additionally, if $f$ is bijective and $f^{-1}$ is also a homomorphism, that is, if a set of vertices in $V_1$ is a hyperedge if and only if their images form a hyperedge of the same size in $H_2$, then it is known as an \emph{isomorphism}. If $H_1=H_2$, then a homomorphism $f:H_1\rightarrow H_2$ is known as an \emph{endomorphism}, and additionally, if $f$ is an isomorphism, it is known as an \emph{automorphism}. A hypergraph $H$ is called a \emph{core} if every endomorphism of $H$ is an automorphism. For an arbitrary hypergraph $H$, a \emph{core of $H$}, denoted $H^\bullet$, is a subhypergraph of $H$ such that $H^\bullet$ is a core and there exists a homomorphism $H \rightarrow H^\bullet$. Figure \ref{example} shows $H^\bullet$, the core of $H$, such that there is a homomorphism $h:H \rightarrow H^\bullet$ defined by $h(a)=h(e)=e, h(b)=b, h(c)=c, h(d)=d$, $h(e)=e$, $h(f)=f$, $h(g)=d$. \\
 
A hypergraph is \emph{vertex-transitive} if its automorphism group acts transitively on its vertices, that is, for any pair of vertices, there is an automorphism taking one vertex to the other. Cores of vertex-transitive graphs have been studied in \cite{Roberson_2013}. A vertex-transitive graph is a core, if its order is relatively prime to its independence number \cite{hahn1997graph}. We verify that this sufficient condition can be extended to hypergraphs as well, using the following lemmas.
\begin{lemma}\cite{Raaphorst2013}\label{vthypercore}
    The core $H^\bullet$ of a vertex-transitive hypergraph $H$ is also vertex-transitive.
\end{lemma}
A set of vertices in a hypergraph $H$ is said to be \emph{strongly independent} if no two belong to the same hyperedge in $H$. The maximum cardinality of such a set is the \emph{strong independence number}, denoted $\alpha_S(H)$. In a graph $G$, these terms correspond to the usual \emph{independent set} and \emph{independence number}, $\alpha(G)$, respectively. For example, in the hypergraph $H$ from Figure~\ref{example}, $\{a,d,g\}$ is a strongly independent set and $\alpha_S(H)=3$. The next result is an extension of the No-Homomorphism lemma from \cite{albertson1985homomorphisms}. 
\begin{lemma}\label{nohomo}
    Let $H_1=(V_1, E_1), H_2=(V_2, E_2)$ be hypergraphs such that $H_1$ is vertex-transitive and there exists a homomorphism $H_2 \rightarrow H_1$. Then, $$ \frac{\alpha_S(H_2)}{|V_2|} \geq \frac{\alpha_S(H_1)}{|V_1|}.$$
\end{lemma}
\begin{proof}
    Let $\mathcal{S}(H_1)$ denote the family of strongly independent sets of size $\alpha_S(H_1)$ in $H_1$. Since $H_1$ is vertex-transitive, for any pair of vertices $u,v\in V$, there exists an automorphism $\phi$ on $H_1$ such that $\phi(u)=v$. Also, under any automorphism $\phi$, strongly independent sets are mapped onto strongly independent sets. 
    So, for any maximum strongly independent set $I_1$ containing $u$, there exists an automorphism $\phi$ on $H$ and a strongly independent set $J_1$ containing $v$ such that $\phi(I_1)=J_1$. Let $I_1, I_2, \ldots , I_l$ be those members of $\mathcal{S}(H_1)$ containing $u$ and  $J_1, J_2, \ldots , J_k$ be those members of $\mathcal{S}(H_1)$ containing $v$. Then $\phi$ being an automorphism takes each $I_i$ to a distinct $J_j$, so $l\leq k$. Similarly $\phi^{-1}$ takes each $J_j$ to a distinct $I_i$, so $l\geq k$ and hence $l=k$. Therefore, every vertex of $H_1$ is in the same number, say $m$, of the members of $\mathcal{S}(H_1)$. Hence, counting the cardinalities of all the members of $\mathcal{S}(H_1)$ in two ways, we have \begin{equation}\alpha_S(H_1)\cdot |\mathcal{S}(H_1)| = |V_1|\cdot m.\end{equation}
    Now, for any homomorphism $\phi ': H_2 \rightarrow H_1$ and for each $I \in \mathcal{S}(H_1)$, the cardinality of the pre-image satisfies: $|\phi'^{-1}(I)| \leq  \alpha_S(H_2)$.
Summing this inequality over all the members of $\mathcal{S}(H)$, we get
\begin{equation}\sum\limits_{I\in \mathcal{S}(H_1)}
|\phi '^{-1}(I)| \leq \alpha_S(H_2) \cdot |\mathcal{S}(H_1)|\end{equation}
However, each $u \in V (H_2)$ contributes exactly $m$ to the sum $\sum\limits_{I\in \mathcal{S}(H_1)}
|\phi '^{-1}(I)|$, since $\phi '(u)$ belongs to exactly $m$ members of $\mathcal{S}(H_1)$. Thus,
\begin{equation}\sum\limits_{I\in \mathcal{S}(H_1)}
|\phi '^{-1}(I)| = m\cdot |V_2|\end{equation}
Combining equations (1), (2) and (3), we get
$$ \alpha_S(H_2)/|V_2| \geq m/|\mathcal{S}(H_1)| = \alpha_S(H_1)/|V_1|.$$
\end{proof}
\begin{theorem}\label{mainref}
    Let $H=(V, E)$ be a vertex-transitive hypergraph. If the order of the hypergraph, $|V|$ and the strong independence number, $\alpha_S(H)$ are relatively prime, then $H$ is a core.
\end{theorem}
\begin{proof}
    Suppose $H=(V, E)$ satisfies the given conditions, but is not a core. Let $H^\bullet=(V^\bullet, E^\bullet)$ be a proper subhypergraph which is the core of $H$. Then, by Lemmas \ref{vthypercore} and \ref{nohomo}, $\frac{\alpha_S(H)}{|V|} = \frac{\alpha_S(H^\bullet)}{|V^\bullet|}$. Since $|V|$ and $\alpha_S(H)$ are relatively prime, $\frac{\alpha_S(H)}{|V|}$ is in the simplest form. This contradicts the assumption that $|V^\bullet|<|V|$. Hence, $H$ is a core.
\end{proof}

\begin{lemma}\cite{raina2026coresuniformalmostuniform3qualitative}\label{VT}
    The hypergraphs $3\text{-}UQI(n, 2)$ and $3\text{-}AUQI(n, 2)$ are vertex-transitive.
\end{lemma}
\section{Strong Independence Number of $3\text{-}AUQI(2n+1,2)$}\label{sec:3}
  In this section, we establish that the strong independence number of the almost-uniform qualitative independence hypergraph, $\alpha_S(3\text{-}AUQI(2n+1,2))=2n+3$, for $n\geq 4$, subsequently giving a lower bound on the strong chromatic number of $3\text{-}QI(2n+1,2)$. Since the vertices of $3\text{-}AUQI(2n+1,2)$ are binary strings, we identify them with subsets of $\{1,2,\ldots,2n+1\}$ of size $n$ in the corresponding partition. We use the following condition to first identify the criteria under which two vertices in $3\text{-}AUQI(2n+1,2)$ can belong to a strongly independent set.

\begin{condition}\cite{raina2026coresuniformalmostuniform3qualitative}\label{earlycondition} A pair of distinct vertices $A$ and $B$ cannot belong to a hyperedge in $3\text{-}QI(n, 2)$ if and only if at least one of the four sets $A \cap B$, $A \cap \overline{B}$, $\overline{A} \cap B$, $\overline{A} \cap \overline{B}$ has cardinality at most one.\end{condition}

\begin{condition}\label{conditionchar} A pair of distinct vertices $A$ and $B$ can belong to a strongly independent set in $3\text{-}AUQI(2n+1, 2)$ if and only if $|A \cap B|=0$ or $1$ or $(n-1)$.\end{condition}
\begin{proof}
    $|A\cap B|=n-1$ if and only if $|\overline{A}\cap B|=1$. The proof then follows from Condition \ref{earlycondition}.
\end{proof}
For notational convenience, we drop the curly braces in the set notation of vertices of $3\text{-}AUQI(9,2)$. For example, we denote the vertex $000011111$ of $3\text{-}AUQI(9,2)$ by $1234$.
\begin{lemma}\label{int1}
     Let $A=a_1a_2\ldots a_n$ be a vertex of $3\text{-}AUQI(2n+1,2)=(V, E)$. Consider a family $\mathcal{F}=\{B\in V: |A \cap B|=1\}$. Then, at most $n$ from $\mathcal{F}$ along with $A$ can form a strongly independent set in $3\text{-}AUQI(2n+1,2)$, when $n>4$ and at most $6$ from $\mathcal{F}$ along with $A$ can form a strongly independent set in $3\text{-}AUQI(2n+1,2)$, when $n=4$.
     
\end{lemma}
\begin{proof}
    Let $\mathcal{S}\subset \mathcal{F}$ be such that it forms a strongly independent set along with $A$. Without loss of generality, let $B=a_1a_{n+1}a_{n+2}\ldots a_{2n-1} \in \mathcal{S}$. For any $B_1\neq B$ in $\mathcal{S}$, such that $|B\cap B_1|=n-1$, we have two cases:\\

    \noindent \textbf{Case 1.} If $a_1 \in B\cap B_1$, then the options for $B_1$ are $a_1L_ia_j$, where $i=1,2,\ldots, n-1$, $j=2n, 2n+1$ and $L_i$ is an $(n-2)$-subset of $a_{n+1}a_{n+2}\ldots a_{n+i-1}a_{n+i+1}\ldots a_{2n-1}$. However by the pigeonhole principle, $|a_1L_ia_j \cap a_1L_ka_l|=n-2$ for $i\neq k, j\neq l$. This implies that the only extremal case is $a_1L_ia_j$ for fixed $j$ and $i=1, 2, \ldots, n-1$. Thus, there are $(n-1)$ choices for $B_1$ and therefore, in this case we have $|\mathcal{S}|\leq n$.\\

    \noindent \textbf{Case 2.} If $a_1 \notin B\cap B_1$, then the options for $B_1$ are $ka_{n+1}a_{n+2}\ldots a_{2n-1}$ for $k=a_2, a_3, \ldots, a_n$. Thus, in this case also, we have $|\mathcal{S}|\leq n$.\\

    Since $|a_1L_ia_j \cap ka_{n+1}a_{n+2}\ldots a_{2n-1}|=n-2$, by Condition \ref{conditionchar}, vertices from Case 1 and Case 2 cannot coexist in $\mathcal{S}$. Note that unless $n=4$, for any $B_1\neq B$ in $\mathcal{S}$, we must have $|B\cap B_1|=n-1$. Thus, $|\mathcal{S}|\leq n$ if $n>4$. \\

    Now suppose that $n=4$, $B=a_1a_5a_6a_7$. Then, in addition to the above cases, there can exist a $B_1\neq B$ in $\mathcal{S}$, such that $|B\cap B_1|=1$. We can deduce that $B\cap B_1\neq a_1$, otherwise enough digits are not available to form $B_1$. Therefore, the possible choices of $B_1$ are $a_la_ma_8a_9$, where $l=2,3,4$, $m=5,6,7$. Assume without loss of generality that $B_1=a_2a_5a_8a_9$. This reduces the remaining choices to either of $a_la_5a_8a_9$ for $l=3,4$ or $a_2a_ma_8a_9$ for $m=6,7$. From Case 1, the remaining options are $a_1a_6a_7a_8$ and $a_1a_6a_7a_9$. From Case 2,  the remaining options are $a_3a_5a_6a_7$ and $a_4a_5a_6a_7$. Thus, if $n=4$, $|\mathcal{S}|\leq 1\, (B) + 3 \text{ (intersecting at one digit with $B$) } + 2\text{ (Case 1 or Case 2) }=6$. For example, $1567, 2589, 3589, 4589, 1678, 1679$ form a strongly independent set along with $1234$ in $3\text{-}AUQI(9,2)$.
\end{proof}
We now consider the families of vertices that intersect at $(n-1)$ digits with a given vertex in $3\text{-}AUQI(2n+1,2)$. In order to find the largest subset of such a family that can form a strongly independent set, we use the existence of Latin rectangles. An $r\times n$ \emph{Latin rectangle} based on $[n]$, is an $r\times n$ matrix such that each entry is one of the integers from $[n]$ and each integer occurs in each row and column at most once.
\begin{lemma}\cite{bollobas1998modern}\label{Latin}
    For $r,n \in \mathbb{N}$, $r\leq n$, an $r\times n$ Latin rectangle always exists.
\end{lemma}
\begin{lemma}\label{intn-1}
    Let $A=a_1a_2\ldots a_n$ be a vertex of $3\text{-}AUQI(2n+1,2)=(V, E)$. Consider a family $\mathcal{F}=\{B\in V: |A \cap B|=n-1\}$. Then, at most $(n+1)$ from $\mathcal{F}$ along with $A$ can form a strongly independent set in $3\text{-}AUQI(2n+1,2)$.
\end{lemma}
\begin{proof}
    There are $\binom{n}{n-1}$ ways to choose the intersecting digits from $A$ and $\binom{n+1}{1}$ ways to choose the non-intersecting digit with respect to $A$ to form a vertex in $\mathcal{F}$, so $|\mathcal{F}|=n(n+1)$. We show that at most $(n+1)$ of these can belong to a strongly independent set with $A$. Let $P_i$ denote the $(n-1)$-subset of $A$ omitting $a_i$, for $i=1,2, \ldots, n$. The table on the left in Figure \ref{n-1} shows the vertices in $\mathcal{F}$, where row $i$ corresponds to $P_i$ and column $j$ corresponds to $a_j$, for $i=1,2,\ldots, n$ and $j=n+1, n+2, \ldots, 2n+1$. \\
    
    \begin{figure}
    \begin{center}
        \begin{tabular}{|c|c|c|c|}
        \hline
           $P_1a_{n+1}$  & $P_1a_{n+2}$&$\ldots$&$P_1a_{2n+1}$ \\\hline
             $P_2a_{n+1}$  & $P_2a_{n+2}$&$\ldots$&$P_2a_{2n+1}$ \\\hline 
             $\ldots$&$\ldots$&$\ldots$&$\ldots$\\\hline
             $P_na_{n+1}$  & $P_na_{n+2}$&$\ldots$&$P_na_{2n+1}$ \\\hline
        \end{tabular}\hspace{1cm}
         \begin{tabular}{|c|c|c|c|}
        \hline
           $P_1\boldsymbol{a_{n+1}}$  & $P_1\boldsymbol{a_{n+2}}$&$\ldots$&$P_1\boldsymbol{a_{2n+1}}$ \\\hline
             $P_2\boldsymbol{a_{n+2}}$  & $P_2\boldsymbol{a_{n+3}}$&$\ldots$&$P_2\boldsymbol{a_{n+1}}$ \\\hline 
             $\ldots$&$\ldots$&$\ldots$&$\ldots$\\\hline
             $P_n\boldsymbol{a_{2n}}$  & $P_n\boldsymbol{a_{2n+1}}$&$\ldots$&$P_n\boldsymbol{a_{2n-1}}$ \\\hline
        \end{tabular}
        \caption{Vertices in $\mathcal{F}$, as given in Lemma \ref{intn-1}}
        \label{n-1}\end{center}
    \end{figure}

    We now claim that the entries in each column of the array in Figure \ref{n-1} can be rearranged in such a way that any two vertices in the same column cannot belong to a strongly independent set. If these $(n^2 + n)$ vertices can be organized into $n$ columns such that any two vertices in the same column intersect at exactly $(n-2)$ digits, they would be part of a hyperedge and consequently, cannot belong to a strongly independent set. This arrangement ensures that, at most one vertex from each column can belong to a strongly independent set containing $A$, resulting in a maximum strongly independent set size of $(n+1)$. By the pigeonhole principle, $|P_i \cap P_j|=n-2$ for $i\neq j$. Therefore, in order to establish the claim, the last digit in each vertex in a column must be distinct. If we label each vertex by its last digit, we are looking for a Latin rectangle of order $n \times (n+1)$ with symbols $a_j, j=n+1, n+2, \ldots, 2n+1$, whose existence is guaranteed by Lemma \ref{Latin}. For instance, see the table on the right in Figure \ref{n-1} for a corresponding Latin rectangle rearrangement. Hence, a strongly independent set can have at most $(n+1)$ vertices from $\mathcal{F}$. \\
\end{proof}
\begin{theorem}\label{main}
    The strong independence number of the almost-uniform qualitative independence hypergraph, $3\text{-}AUQI(2n+1,2)=(V, E)$, is given by, $\alpha_S(3\text{-}AUQI(2n+1,2))=2n+3$, for $n\geq 4$.
\end{theorem}
\begin{proof}
    Let $A=a_1a_2\ldots a_n \in V$ and $\overline{A}=a_{n+1}a_{n+2}\ldots a_{2n+1}$. First, we construct a strongly independent set of size $(2n+3)$ from $V$. Consider the families of subsets of $V$, $\mathcal{X}=\{a_1a_2\ldots a_{n-1}a_i: a_i \in \overline{A}\}$ and $\mathcal{Y}$, the set of all $n$-subsets of $\overline{A}$. By definition, all the vertices in $\mathcal{X}$ have $a_1a_2\ldots a_{n-1}$ in common and therefore, any pair of vertices in $\mathcal{X}$ intersect at $(n-1)$ digits. Also, by the pigeonhole principle, any pair of vertices in $\mathcal{Y}$ intersects at $(n-1)$ digits. Let $M\in \mathcal{X}$ and $N\in \mathcal{Y}$ be two vertices. Then, $|A\cap M|=n-1$, $|A\cap N|=0$ and $|M \cap N|\leq 1$. Since $|\mathcal{X}|=|\overline{A}|=n+1$ and $|\mathcal{Y}|=\binom{n+1}{n}=n+1$, by Condition \ref{conditionchar}, $\{A\}\cup \mathcal{X}\cup \mathcal{Y}$ is a strongly independent set of size $(2n+3)$. Thus, $\alpha_S(3\text{-}AUQI(2n+1,2))\geq 2n+3$. \\

    For the reverse inequality, we consider three sub-families of $V$ as suggested by Condition \ref{conditionchar}, namely $\mathcal{F}_1=\{B\in V: |A \cap B|=0\}$, $\mathcal{F}_2=\{B\in V: |A \cap B|=1\}$ and $\mathcal{F}_3=\{B\in V: |A \cap B|=n-1\}$. By the pigeonhole principle, the vertices in $\mathcal{F}_1$, which is the set of all $n$-subsets of $\overline{A}$, intersect pairwise at $(n-1)$ digits. Hence, $|\mathcal{F}_1|=\binom{n+1}{n}=n+1$ and $\mathcal{F}_1$ forms a strongly independent set with $A$. From Lemma \ref{int1}, at most $n$ from $\mathcal{F}_2$ can form a strongly independent set with $A$, if $n>4$ and at most $6$ can form a strongly independent set with $A$ if $n=4$. From Lemma \ref{intn-1}, at most $(n+1)$ from $\mathcal{F}_3$ can form a strongly independent set with $A$.\\

    Let $\mathcal{S}$ be a strongly independent set containing $A$. If $\mathcal{F}_2$ contributes no vertex to $\mathcal{S}$, then $|\mathcal{S}|\leq 1\, (A) + (n+1)\, (\mathcal{F}_1) + (n+1)\, (\mathcal{F}_3)= 2n+3$. Hence, without loss of generality, in the rest of the proof we assume that at least one vertex, say $B=a_1a_{n+1}a_{n+2}\ldots a_{2n-1}$ from $\mathcal{F}_2$, is in $\mathcal{S}$. Since $B$ intersects with any vertex from $\mathcal{F}_1$ at at least $(n-2)$ digits, by Condition \ref{conditionchar}, the only possible choices from $\mathcal{F}_1$ that can be in $\mathcal{S}$ are the ones that intersect at $(n-1)$ digits with $B$, namely, $a_{n+1}a_{n+2}\ldots a_{2n-1}a_i$, for $i=2n, 2n+1$. By Condition \ref{conditionchar}, the options from $\mathcal{F}_3$ given in Figure \ref{n-1} of Lemma \ref{intn-1} reduce to $P_1a_k$ for $k=n+1, n+2, \ldots, 2n+1$ and $P_ja_{2n}, P_ja_{2n+1}$ for $j=2,3, \ldots, n$. Note that $|P_ia_k\cap P_xa_y|=n-2$ for $i\neq x, k\neq y$. So, in the extremal case, $\mathcal{F}_3$ contributes $(n+1)$, only when we take $P_1a_k$ for all $k$. The next-best extremal case here is of $n$ vertices, $P_ja_k$ for fixed $k$, where $k=2n$ or $2n+1$ and $j=1,2,\ldots, n$. If none other than $B$ are chosen from $\mathcal{F}_2$, then $$|\mathcal{S}|\leq 1\, (A) + 1\, (B) + 2\, (\mathcal{F}_1) + (n+1)\, (\mathcal{F}_3)=n+5 < 2n+3,\text{ when }n\geq 3.$$

    Suppose now that in addition to $B$, at least one vertex is chosen from $\mathcal{F}_2$, say $B_1$. If $n>4$, then $B_1$ can be from Case 1 or 2, given in the proof of Lemma \ref{int1}. If $B_1$ is from Case 1 of Lemma \ref{int1}, then only one option remains from $\mathcal{F}_1$ depending on which digit out of $a_{2n}$ and $a_{2n+1}$ is in $B_1$. For instance, if $B_1=a_1a_{n+2}a_{n+3}\ldots a_{2n-1}a_{2n}$, then $C=a_{n+1}a_{n+2}\ldots a_{2n-1}a_{2n+1}$ is eliminated from $\mathcal{F}_1$ since $|B_1\cap C|=n-2$. Thus, $$|\mathcal{S}|\leq 1\, (A) + 1\, (\mathcal{F}_1) + n\, (\mathcal{F}_2) + (n+1)\, (\mathcal{F}_3)=2n+3.$$
    If $B_1$ is from Case 2 of Lemma \ref{int1}, then the extremal case of $P_1a_k$ for all $k=n+1, n+2, \ldots, 2n+1$ from $\mathcal{F}_3$ reduces to just two vertices $P_1a_{2n}$ and $P_1a_{2n+1}$. Therefore, the extremal case now is $P_ja_k$ for fixed $k$ and $j=1,2,\ldots, n$. Thus, $$|\mathcal{S}|\leq 1\, (A) + 2\, (\mathcal{F}_1) + n\, (\mathcal{F}_2) + n\, (\mathcal{F}_3)=2n+3.$$ Finally, let $n=4$ and $B=a_1a_5a_6a_7$. If all the vertices chosen from $\mathcal{F}_2$ are such that they intersect at three digits with $B$, then the proof follows from the previous statements using Cases 1 and 2 from Lemma \ref{int1}. If not, then $B_1$ is chosen from $\mathcal{F}_2$, such that $|B_1\cap B|=1$, say $B_1=a_2a_5a_8a_9$ (without loss of generality). In this case, Lemma \ref{int1} implies that $\mathcal{F}_3$ can contribute up to $6$ vertices to $\mathcal{S}$. However, the choice of $B_1$ reduces the valid choices in $\mathcal{F}_3$ to either $P_2a_8, P_2a_9$ or $P_1a_6, P_1a_7$. Therefore, in this case as well, $$|\mathcal{S}|\leq 1\, (A) + 2\, (\mathcal{F}_1) + 6\, (\mathcal{F}_2) + 2\, (\mathcal{F}_3)=11=2n+3, \text{ for }n=4.$$
    Hence, $\alpha_S(3\text{-}AUQI(2n+1,2))=2n+3$, for $n\geq 4$.
    
\end{proof}
The following is an application of the above result.
\begin{corollary}
    The strong chromatic number of the qualitative independence graph, $3\text{-}QI(2n+1,2)$, is given by, $$\chi_S(3\text{-}QI(2n+1,2))\geq \frac{1}{2n+3}\binom{2n+1}{n}.$$
\end{corollary}
\begin{proof}
    Since $3\text{-}AUQI(2n+1,2)$ is a subgraph of $3\text{-}QI(2n+1,2)$, $\chi_S(3\text{-}QI(2n+1,2))\geq\chi_S(3\text{-}AUQI(2n+1,2))$. The result then follows from Theorem \ref{main}, since $\chi_S(H)\geq \frac{|V|}{\alpha_S(H)}$ for a hypergraph $H=(V, E)$.
\end{proof}
This improves the lower bounds of the strong chromatic number of $3\text{-}QI(9,2)$ given in \cite{raina} as $\chi_S(3\text{-}QI(9,2))\geq 11$ to $\chi_S(3\text{-}QI(9,2))\geq 12$ and of $3\text{-}QI(11,2)$ given in \cite{raina2026strongcolouringqualitativeindependence} $\chi_S(3\text{-}QI(11,2))\geq 31$ to $\chi_S(3\text{-}QI(11,2))\geq 36$.
\section{$3\text{-}AUQI(2n+1,2)$ is a Core for Infinitely Many $n$}\label{sec:4}
In this section, we give a sufficient condition for the almost-uniform qualitative independence hypergraph $3\text{-}AUQI(2n+1,2)$ to be a core for $n\geq 4$ using its strong independence number, which we determined in Section \ref{sec:3}. We then go on to establish that there are cases in which the sufficient condition actually holds by giving an instance of an infinite family satisfying the given condition.
\begin{theorem}\label{sufficiency}
    For $n\geq 4$, if $\binom{2n+1}{n}$ and $(2n+3)$ are relatively prime, then the almost-uniform qualitative independence hypergraph $3\text{-}AUQI(2n+1,2)$ is a core.
\end{theorem}
\begin{proof}
    This follows immediately from Theorem \ref{mainref} since by Lemma \ref{VT}, $3\text{-}AUQI(2n+1,2)$ is vertex-transitive and by Theorem \ref{main}, $\alpha_S(3\text{-}AUQI(2n+1,2))=2n+3$.
\end{proof}
\begin{corollary}
    The almost-uniform qualitative independence hypergraph $3\text{-}AUQI(p-2,2)$ is a core for any prime $p\geq 11$.
\end{corollary}
\begin{proof}
    Consider $p$ and $\dbinom{p-2}{\frac{p-3}{2}}$, for a prime $p\geq 11$. If they are not relatively prime, then $p$ must divide $\dbinom{p-2}{\frac{p-3}{2}}$. This implies that $p$ must divide $(p-2)!$, which is not true by Euclid's lemma, since $p$ is a prime. Thus, $p$ and $\dbinom{p-2}{\frac{p-3}{2}}$ are relatively prime and therefore, by Theorem \ref{sufficiency}, $3\text{-}AUQI(p-2,2)$ is a core.
\end{proof}
\begin{corollary}
    The qualitative independence hypergraph, $3\text{-}QI(9,2)$, is a core.
\end{corollary}

Since there is an infinite collection of primes $p\geq 11$, this establishes that $3\text{-}AUQI(p-2,2)$ is a core, for infinitely many $p$. Consequently, we also prove that an infinite family of merged Johnson graphs are cores. Let $n, m, i$ be fixed positive integers with $n\geq m \geq i$. Define \emph{generalized Johnson graph}, also known as \emph{uniform subset graph} or \emph{distance graph}, denoted $J(n, m, i)$, as follows: The vertices of $J(n, m, i)$ are the $m$-subsets of $[n]$, where two subsets are adjacent if their intersection has size $i$. If $I$ is any subset of $[m]$, the \emph{merged Johnson graph}, $J (n,m)_I$, is defined to be the graph with the vertex set, all the $m$-subsets of $[n]$, with vertices $M$ and $M'$ adjacent if and only if $|M \cap M'| = i$ for some $i \in I$. The $2$-section of $3\text{-}AUQI(n,2)$ is the merged Johnson graph, $J(n, \frac{n-1}{2})_{I}$, that is, $[3\text{-}AUQI(n,2)]_2 = J(n, \frac{n-1}{2})_{I}$, where $n$ is odd and $I=\{2,3, \ldots, \frac{n-1}{2}-2\}$ (Lemma 15, \cite{raina2026coresuniformalmostuniform3qualitative}). Since the vertex set of the merged Johnson graphs $J(2n+1, n)_{I}$ for $I=\{2,3, \ldots, \frac{n-1}{2}-2\}$, and almost-uniform qualitative independence graphs $3\text{-}AUQI(2n+1,2)$, $n\geq 4$ are the same, and $\alpha_S(H)=\alpha([H]_2)$, the proofs of the results in Sections \ref{sec:3} and \ref{sec:4} establish the following results. 
\begin{corollary}
    The independence number of the merged Johnson graph $J(2n+1, n)_{I}$ for $I=\{2,3, \ldots, \frac{n-1}{2}-2\}$, $\alpha(J(2n+1, n)_{I})=2n+3$.
\end{corollary}
\begin{corollary}
    The merged Johnson graph $J(2n+1, n)_{I}$ for $I=\{2,3, \ldots, \frac{n-1}{2}-2\}$ is a core whenever $(2n+3)$ is a prime greater than $7$.
\end{corollary}
\section{Future directions}\label{sec:5}
We have made considerable progress towards determining the core of the almost-uniform qualitative independence hypergraph $3\text{-}AUQI(2n+1,2)$, by demonstrating an infinite family of numbers for which it is a core. To the best of our knowledge, post Meagher and Steven's work \cite{meagher2005covering}, this is the first instance that establishes that an infinite sub-class of almost-uniform $t$-QI hypergraphs are cores. However, there are unresolved cases for which the problem of determining whether $3\text{-}AUQI(2n+1,2)$ is a core or not is open, where $(2n+3)$ is not a prime. It is also an open question to determine whether the uniform qualitative independence hypergraph $3\text{-}UQI(2n,2)$ is a core. Only one case, where $n=4$ has been settled \cite{raina2026coresuniformalmostuniform3qualitative}. A related pair of questions are whether $3\text{-}UQI(n, 2)$ is the core of $3\text{-}QI(n, 2)$ for even $n>8$, and whether $3\text{-}AUQI(n, 2)$ is the core of $3\text{-}QI(n, 2)$ for odd $n>9$. Answering this would address an even more important problem: Does the existence of a covering array on a hypergraph with given parameters guarantee the existence of a balanced covering array on the hypergraph with the same parameters?\\

\noindent\textbf{Acknowledgements }The first author thanks the Department of Science and Technology (DST) for financial support under the DST-INSPIRE Senior Research Fellow scheme with sanction no. DST/INSPIRE Fellowship/2021/IF210377.

\bibliographystyle{plainurl}

\bibliography{ref}

\end{document}